\documentclass[11pt]{amsart}
\usepackage{amsmath,amssymb,amsthm,mathtools}
\usepackage[margin=1.15in]{geometry}
\usepackage{enumitem,booktabs}
\usepackage[dvipsnames]{xcolor}
\usepackage[colorlinks=true,linkcolor=blue!50!black,citecolor=blue!50!black,urlcolor=blue!50!black]{hyperref}

\newtheorem{theorem}{Theorem}[section]
\newtheorem{lemma}[theorem]{Lemma}
\newtheorem{proposition}[theorem]{Proposition}
\newtheorem{corollary}[theorem]{Corollary}
\newtheorem{conjecture}[theorem]{Conjecture}
\theoremstyle{definition}
\newtheorem{definition}[theorem]{Definition}

\theoremstyle{remark}
\newtheorem{remark}[theorem]{Remark}

\newcommand{\R}{\mathbb{R}}
\newcommand{\Q}{\mathbb{Q}}
\newcommand{\N}{\mathbb{N}}
\newcommand{\Z}{\mathbb{Z}}
\newcommand{\X}{\mathcal{X}}
\newcommand{\RI}{\mathcal{R}}
\newcommand{\cP}{\mathcal{P}}
\newcommand{\cH}{\mathcal{H}}
\newcommand{\one}{\mathbf{1}}

\newcommand{\sd}{\mathbin{\triangle}}

\title{Chord theorems on graphs (revised and extended)}
\author{Mohammad Javaheri}
\address{Department of Mathematics and Statistics, Siena University\\
515 Loudon Road, Loudonville, NY 12211, USA}
\email{mjavaheri@siena.edu}

\begin{document}
\begin{abstract}
The Horizontal Chord Theorem states that if a continuous curve connects
points $A$ and $B$ in the plane, then for every positive integer $k$ there
are points $C$ and $D$ on the curve with $\overrightarrow{AB}=k\,\overrightarrow{CD}$.
We reformulate the one-dimensional chord problem on finite metric graphs:
the \emph{chord set} $\cH(G)$ consists of those $r$ such that every function
of mean zero on $G$ integrates to zero over some closed connected subset of
$G$ of measure $r$. We prove that the space of connected subsets of fixed
measure is compact and path-connected, that $\cH(G)$ always contains the
closure of a combinatorially defined \emph{partition set} $\cP(G)$, and that
the two sets agree at every rational number; as a consequence, membership of
a rational number in $\cH(G)$ is decided by a finite linear program, and the
conjecture $\cH(G)=\overline{\cP(G)}$ reduces to a statement about irrational
points. We give a covering lemma that yields $\cH(G)=[0,|E|]$ for Eulerian
graphs, $[0,1]\subseteq\cH(G)$ when $G$ has no vertex of degree one, and
$[0,g(G)]\subseteq\cH(G)$ when $G$ is bridgeless with girth $g(G)$. Conversely,
a bridge forces a gap in $\cH(G)$ below $|E|$. Computations on small graphs
support the conjecture that $\cH(G)=[0,|E|]$ exactly when $G$ is bridgeless.
\end{abstract}
\maketitle

\section{Introduction}

For a function $F:\R\to\R$ the \emph{chord set} is
$H(F)=\{t\in\R:\ F(x+t)=F(x)\text{ for some }x\}$. Lévy's Universal Chord
Theorem \cite{Levy} says that if $F$ is continuous on $[0,1]$ with
$F(0)=F(1)$, then $1/k\in H(F)$ for every $k\in\N$, and that for every
$r\in(0,1)$ that is not of the form $1/k$ there is such an $F$ with
$r\notin H(F)$. Hopf \cite{Hopf} extended this to planar curves: the set of
chord vectors of a continuous curve from $A$ to $B$ that are parallel to
$\overrightarrow{AB}$ always contains $\frac1k\overrightarrow{AB}$, and its
complement in $(0,1)\overrightarrow{AB}$ is closed under addition. See
Oxtoby \cite{Oxtoby} and Totik \cite{Totik} for surveys and applications.

Writing $F(x)=\int_0^x f$, the statement $r\in H(F)$ for all $F$ becomes:
\emph{every $f$ with $\int_0^1 f=0$ integrates to zero over some interval of
length $r$.} Intervals are exactly the closed connected subsets of the graph
with one edge, and this formulation makes sense on any finite graph. The
purpose of this note is to study the resulting chord sets.

\subsection*{Main results}
Let $G=(V,E)$ be a finite connected graph (loops and multiple edges allowed),
viewed as a metric measure space in which every edge is a copy of $[0,1]$.
Let $\X_r(G)$ be the family of closed connected subsets of $G$ of measure
$r$. Definitions of the chord set $\cH(G)$ and the partition set $\cP(G)$ are
given in Section~\ref{sec:sets}.
\begin{enumerate}[label=(\roman*)]
\item $\X_r(G)$ is compact and path-connected (Theorem~\ref{thm:pathconn}).
\item $\overline{\cP(G)}\subseteq\cH(G)$ and $\cH(G)$ is closed
  (Theorem~\ref{thm:PH}, Proposition~\ref{prop:closed}).
\item $\cP(G)\subseteq\Q$ and $\cH(G)\cap\Q=\cP(G)$; for $r=p/q$ in lowest
  terms, membership is decided by a finite linear program over
  $\frac1q$-grid-aligned sets (Theorem~\ref{thm:rational}).
\item $r\in\cH(G)$ if and only if the constant function $\one$ is an
  $L^1$-limit of nonnegative combinations of indicators of sets in $\X_r(G)$
  (Theorem~\ref{thm:dual}). Hence the conjecture
  $\cH(G)=\overline{\cP(G)}$ is
  equivalent to: \emph{every point of $\cH(G)$ is a limit of rational points
  of $\cH(G)$.}
\item A covering lemma (Lemma~\ref{lem:covering}) gives $\cH(G)=[0,|E|]$ for
  Eulerian $G$, $[0,1]\subseteq\cH(G)$ if $G$ has minimum degree at least
  $2$, and $[0,g]\subseteq\cH(G)$ if $G$ is bridgeless of girth $g$
  (Theorems~\ref{thm:mindeg} and \ref{thm:girth}).
\item If $G$ has a bridge then $\cH(G)$ misses an interval
  $(\mu,|E|)$ with $\mu<|E|$ given explicitly (Proposition~\ref{prop:bridge}).
  We conjecture that $\cH(G)=[0,|E|]$ if and only if $G$ is bridgeless
  (Conjecture~\ref{conj:bridgeless}); Section~\ref{sec:comp} reports exact
  computations supporting this.
\end{enumerate}
Section~\ref{sec:necklace} returns to the motivating necklace problem and to
a combinatorial game, where we also correct the graph version of the game.

\section{Connected subsets of a metric graph}\label{sec:prelim}

Each edge $e$ of $G$ is identified with $[0,1]$ (for a loop, both ends are
identified with the same vertex). This makes $G$ a compact length space; we
write $d_G$ for its path metric and $\mu$ for the measure that is Lebesgue
measure on each edge, so $\mu(G)=|E|$. We write $\RI(G)$ for the set of
functions $f:G\to\R$ whose restriction to each edge is Riemann integrable,
and $I_f(A)=\int_A f\,d\mu$. Every $f\in\RI(G)$ is bounded; put
$\|f\|_\infty=\sup|f|$.

Let $\X(G)$ be the set of nonempty closed connected subsets of $G$ and
$\X_r(G)=\{A\in\X(G):\mu(A)=r\}$ for $0<r<|E|$.

\begin{lemma}[Structure]\label{lem:structure}
Let $A\in\X(G)$ and $S=A\cap V$.
\begin{enumerate}[label=(\alph*)]
\item If $S=\emptyset$, then $A$ is a point or a closed interval in the
  interior of a single edge.
\item If $S\neq\emptyset$, then for every edge $e=[u,v]$ the set $A\cap e$
  is: empty if $u,v\notin S$; an initial segment $[0,a]$ with $0\le a<1$ if
  $u\in S$, $v\notin S$; and $[0,a]\cup[b,1]$ with $0\le a\le b\le 1$ if
  $u,v\in S$ (this includes $A\cap e=e$). Moreover $S$ is connected through
  the edges contained in $A$.
\end{enumerate}
In particular, if $\mu(A)>0$ then $A$ is the closure of the union of its
nondegenerate pieces.
\end{lemma}

\begin{proof}
Let $J$ be a connected component of $A\cap e^\circ$, where $e^\circ$ is the
interior of $e$. If $\overline J$ contains no point of $S$, then $J$ is
closed in $G$, and it is open in $A$ because $J=A\cap U$ for a small open
neighborhood $U$ of $\overline J$ in $e^\circ$. Since $A$ is connected,
$A=J$, which is case (a). Otherwise every component of $A\cap e^\circ$ has
an endpoint of $e$ in $S$ in its closure, which gives the list in (b). The
pieces $[0,a]$ attached at a single vertex do not join distinct vertices,
so the vertices of $S$ are joined inside $A$ only through full edges.
For the last claim, a vertex $v\in S$ that is not in the closure of a
nondegenerate piece has only degenerate pieces at it and no full edge, so
$A=\{v\}$ by connectedness.
\end{proof}

For $A,B\in\X_r(G)$ put
\[
\delta(A,B)=\mu(A\sd B)=2r-2\mu(A\cap B).
\]
By the last assertion of Lemma~\ref{lem:structure}, $\delta(A,B)=0$ implies
$A=B$, so $\delta$ is a metric on $\X_r(G)$.

\begin{lemma}\label{lem:Icont}
For $f\in\RI(G)$ and $A,B\in\X_r(G)$, $|I_f(A)-I_f(B)|\le\|f\|_\infty\,\delta(A,B)$.
\end{lemma}

\begin{proof}
$|I_f(A)-I_f(B)|=\bigl|\int_{A\setminus B}f-\int_{B\setminus A}f\bigr|\le\int_{A\sd B}|f|$.
\end{proof}

\begin{lemma}\label{lem:compact}
$(\X_r(G),\delta)$ is compact.
\end{lemma}

\begin{proof}
By Lemma~\ref{lem:structure}, a set $A\in\X(G)$ is determined by finitely
many combinatorial data (the set $S$, the edges contained in $A$, the
pattern of the remaining pieces) and finitely many real parameters in
$[0,1]$ (the endpoints of the pieces). Given a sequence in $\X_r(G)$, pass
to a subsequence with constant combinatorial data and convergent
parameters. The set $A$ described by the limit parameters is closed and
connected (a piece $[0,a_j]$ may converge to a full edge or to a vertex,
which preserves connectedness), it has measure $r$ because $\mu$ is the
sum of the piece lengths, and $\delta(A_j,A)\to0$ because $\delta$ is
bounded by the sum of the changes in the parameters.
\end{proof}

\begin{lemma}[Retraction and growth]\label{lem:retract}
Let $C\in\X(G)$ with $\mu(C)=c$.
\begin{enumerate}[label=(\alph*)]
\item For every $x\in C$ there is a family $\{C\langle a\rangle\}_{0\le a\le c}$
  in $\X(G)$ with $C\langle 0\rangle=\{x\}$, $C\langle c\rangle=C$,
  $\mu(C\langle a\rangle)=a$, and $C\langle a\rangle\subseteq C\langle a'\rangle$
  for $a\le a'$. In particular $\mu(C\langle a\rangle\sd C\langle a'\rangle)=|a-a'|$.
\item There is a family $\{C^{(a)}\}_{c\le a\le|E|}$ in $\X(G)$ with
  $C^{(c)}=C$, $C^{(|E|)}=G$, $\mu(C^{(a)})=a$, increasing in $a$.
\end{enumerate}
\end{lemma}

\begin{proof}
(a) Let $\rho(y)$ be the infimum of the lengths of paths in $C$ from $x$
to $y$. By Lemma~\ref{lem:structure}, $C$ is a finite union of segments, so
$\rho$ is continuous and piecewise linear with slopes $\pm1$; in particular
each level set $\{\rho=s\}$ is finite. The sets $C^s=\{\rho\le s\}$ are
closed, connected (each point is joined to $x$ by a shortest path inside
$C^s$), and increasing in $s$. The function $m(s)=\mu(C^s)$ is continuous
because level sets are null, and strictly increasing on $[0,\max\rho]$:
if $s<s'\le\max\rho$, a shortest path from $x$ to a point at level $s'$
contributes length $s'-s$ to $C^{s'}\setminus C^s$. Since $m(0)=0$ and
$m(\max\rho)=c$, we may set $C\langle a\rangle=C^{m^{-1}(a)}$.

(b) Use the same construction with $\rho(y)=d_G(y,C)$ and $C^s=\{\rho\le s\}$,
which is connected because every point is joined to $C$ by a shortest path.
\end{proof}

\section{Path-connectedness}\label{sec:pathconn}

\begin{theorem}\label{thm:pathconn}
For $0<r<|E|$, the space $(\X_r(G),\delta)$ is compact and path-connected.
\end{theorem}

\begin{proof}
Compactness is Lemma~\ref{lem:compact}. Let $A,B\in\X_r(G)$.

\emph{Case 1: $A\cap B\neq\emptyset$.} Fix $x\in A\cap B$ and let
$A\langle a\rangle$, $B\langle b\rangle$ be the retractions of
Lemma~\ref{lem:retract}(a) towards $x$. For $(a,b)\in[0,r]^2$ the set
$\Phi(a,b)=A\langle a\rangle\cup B\langle b\rangle$ is closed and connected
(both pieces contain $x$), and
\[
\delta\bigl(\Phi(a,b),\Phi(a',b')\bigr)\le|a-a'|+|b-b'| .
\]
Let $g(a,b)=\mu(\Phi(a,b))$. Then $g$ is continuous, nondecreasing in each
variable, $g(a,0)=a$, $g(0,b)=b$ and $g(a,r)\ge r$. Put
$\beta(a)=\min\{b\in[0,r]:g(a,b)\ge r\}$. If $\beta(a)>0$ then
$g(a,b)<r$ for $b<\beta(a)$, so $g(a,\beta(a))=r$ by continuity; if
$\beta(a)=0$ then $a=g(a,0)\ge r$, so $a=r$ and again $g=r$. Monotonicity
of $g$ in $a$ shows that $\beta$ is nonincreasing, with $\beta(r)=0$ and
$\beta(0)=r$. Let
\[
\Gamma=\{(a,b)\in[0,r]^2:\ \beta(a^+)\le b\le\beta(a^-)\},
\]
the graph of $\beta$ with its jumps filled in by vertical segments (with
$\beta(0^-)=r$, $\beta(r^+)=0$). Passing to limits shows $g=r$ at the two
ends of each vertical segment, and then on the whole segment by
monotonicity in $b$; so $g\equiv r$ on $\Gamma$. The function $(a,b)\mapsto a-b$
is continuous and injective on the compact connected set $\Gamma$, so
$\Gamma$ is an arc from $(0,r)$ to $(r,0)$. Along $\Gamma$, $\Phi$ is a
continuous path in $\X_r(G)$ from $\Phi(0,r)=B$ to $\Phi(r,0)=A$.

\emph{Case 2: $A\cap B=\emptyset$.} Let $\gamma:[0,\ell]\to G$ be a shortest
path from $A$ to $B$, parametrized by arc length, so that $\gamma$ is
injective, $\gamma(0)\in A$, $\gamma(\ell)\in B$, and $\gamma((0,\ell))$
misses $A\cup B$. Let $A\langle\cdot\rangle$ be the retraction of $A$ towards
$\gamma(0)$. For $0\le t\le\min(r,\ell)$ the sets
$U_t=A\langle r-t\rangle\cup\gamma([0,t])$ lie in $\X_r(G)$ and depend
continuously on $t$. If $\ell\le r$, then $U_\ell$ meets $B$ and we conclude
by Case~1. If $\ell>r$, then $U_r=\gamma([0,r])$, and
$V_t=\gamma([t,t+r])$, $0\le t\le\ell-r$, is a continuous path in
$\X_r(G)$ from $U_r$ to a set meeting $B$; again Case~1 applies.
\end{proof}

\section{Partition sets and chord sets}\label{sec:sets}

\begin{definition}\label{def:sets}
The \emph{partition set} $\cP(G)$ is the set of $r\in(0,|E|)$ for which
there are $n\in\N$ and a finite family $U_1,\dots,U_k\in\X_r(G)$
(repetitions allowed) such that every point of $G$ outside a finite set
lies in exactly $n$ of the $U_\alpha$; equivalently
$\sum_\alpha\one_{U_\alpha}=n$ almost everywhere.

The \emph{chord set} $\cH(G)$ is the set of $r\in[0,|E|]$ such that for every
$f\in\RI(G)$ with $\int_G f=0$ there is $U$ with $\mu(U)=r$, $U$ closed and
connected, and $I_f(U)=0$.
\end{definition}

Trivially $0,|E|\in\cH(G)$ (take a point, or $G$). The following
reformulation is used repeatedly.

\begin{lemma}[Sign principle]\label{lem:sign}
Let $0<r<|E|$ and $f\in\RI(G)$. If there are $U,V\in\X_r(G)$ with
$I_f(U)\ge0\ge I_f(V)$, then $I_f(W)=0$ for some $W\in\X_r(G)$.
Consequently $r\notin\cH(G)$ if and only if there is $f\in\RI(G)$ with
$\int_Gf=0$ and $I_f>0$ on $\X_r(G)$; by compactness, then $I_f\ge\varepsilon>0$
on $\X_r(G)$.
\end{lemma}

\begin{proof}
$I_f$ is continuous (Lemma~\ref{lem:Icont}) on the path-connected,
compact space $\X_r(G)$ (Theorem~\ref{thm:pathconn}); apply the
intermediate value theorem along a path from $U$ to $V$. Replacing $f$ by
$-f$ if necessary gives the second statement.
\end{proof}

\begin{proposition}\label{prop:closed}
$\cH(G)$ is closed.
\end{proposition}

\begin{proof}
Let $r_j\to r$ with $r_j\in\cH(G)$, and $0<r<|E|$ (the endpoints are
trivial). Let $\int_Gf=0$ and choose $U_j\in\X_{r_j}(G)$ with $I_f(U_j)=0$.
By Lemma~\ref{lem:retract} there are $U_j'\in\X_r(G)$ with
$\mu(U_j\sd U_j')=|r_j-r|$, so $|I_f(U_j')|\le\|f\|_\infty|r_j-r|\to0$. By
compactness a subsequence of $U'_j$ converges to some $U\in\X_r(G)$ with
$I_f(U)=0$.
\end{proof}

\begin{theorem}\label{thm:PH}
For every finite connected graph $G$, $\overline{\cP(G)}\subseteq\cH(G)$.
\end{theorem}

\begin{proof}
By Proposition~\ref{prop:closed} it suffices to show $\cP(G)\subseteq\cH(G)$.
If $\sum_\alpha\one_{U_\alpha}=n$ a.e.\ and $\int_Gf=0$, then
$\sum_\alpha I_f(U_\alpha)=n\int_G f=0$, so some $I_f(U_\alpha)\ge0$ and some
$I_f(U_\beta)\le0$; apply Lemma~\ref{lem:sign}.
\end{proof}

\subsection{Rational points}
Integrating the identity $\sum_{\alpha=1}^k\one_{U_\alpha}=n$ gives
$kr=n|E|$, so
\begin{equation}\label{eq:Prational}
\cP(G)\subseteq\Q .
\end{equation}

Fix $q\in\N$ and divide every edge into $q$ cells of length $1/q$. A
\emph{$q$-step function} is a function constant on the interior of every
cell, and a set $U\in\X(G)$ is \emph{$q$-aligned} if all endpoints of its
pieces (Lemma~\ref{lem:structure}) lie in $\frac1q\Z$. For $r\in\frac1q\Z$
the set $\X_r^{(q)}(G)$ of $q$-aligned members of $\X_r(G)$ is finite.

\begin{lemma}[Grid lemma]\label{lem:grid}
Let $r\in\frac1q\Z\cap(0,|E|)$ and let $f$ be a $q$-step function. Then for
each $U\in\X_r(G)$, $I_f(U)$ is a convex combination of values $I_f(U_i)$
with $U_i\in\X_r^{(q)}(G)$. In particular
$\min_{\X_r(G)}I_f=\min_{\X_r^{(q)}(G)}I_f$.
\end{lemma}

\begin{proof}
If $U\cap V=\emptyset$, then $U=[x,x+r]$ inside one edge. Since
$r\in\frac1q\Z$, both ends of $U$ cross grid points simultaneously, so
$x\mapsto I_f([x,x+r])$ is affine between consecutive grid values of $x$.

Otherwise, fix the combinatorial data of $U$ (Lemma~\ref{lem:structure}):
the vertex set $S$, the set $F$ of edges contained in $U$, and the list of
remaining pieces $[0,t_i]$ (``stubs''), each attached at a vertex of $S$.
Every vector $t$ in
\[
Q=\Bigl\{t\in[0,1]^m:\ \textstyle\sum_i t_i=r-|F|,\ \ t_i+t_j\le1
\text{ whenever $i,j$ are the two ends of one edge}\Bigr\}
\]
defines a set $U(t)\in\X_r(G)$ (it is connected because $S$ is connected
through $F$). On $Q$, $I_f(U(t))=c+\sum_i\varphi_i(t_i)$ where each
$\varphi_i$ is piecewise affine with breakpoints in $\frac1q\Z$. Hence $I_f$
is affine on each cell $Q_a=Q\cap\prod_i[a_i/q,(a_i+1)/q]$. Order the
coordinates so that the two ends of the same edge are adjacent. Then every
constraint row of $Q_a$ (the coordinate rows, the pair rows and the
all-ones row) has the consecutive-ones property, so the constraint matrix
is an interval matrix and therefore totally unimodular \cite{Schrijver}.
All right-hand sides lie in $\frac1q\Z$, so the vertices of $Q_a$ lie in
$\frac1q\Z^m$, i.e.\ they correspond to $q$-aligned sets. Write $t$ as a
convex combination of the vertices of its cell and use that $I_f$ is affine
there.
\end{proof}

\begin{theorem}\label{thm:rational}
Let $r=p/q\in(0,|E|)$ with $\gcd(p,q)=1$. The following are equivalent:
\begin{enumerate}[label=(\alph*)]
\item $r\in\cH(G)$;
\item $r\in\cP(G)$;
\item $\one$ is a nonnegative combination of the indicators of the finitely
  many sets in $\X_r^{(q)}(G)$.
\end{enumerate}
Consequently $\cH(G)\cap\Q=\cP(G)$, and membership of a rational number in
$\cH(G)$ is decided by a finite linear program.
\end{theorem}

\begin{proof}
(c)$\Rightarrow$(b): the system has rational data, so if it is feasible it
has a rational solution; clearing denominators gives a family with
repetitions covering $G$ exactly $n$ times. (b)$\Rightarrow$(a) is
Theorem~\ref{thm:PH}.

(a)$\Rightarrow$(c): suppose (c) fails. Index the $q|E|$ cells and view each
$\one_U$, $U\in\X_r^{(q)}(G)$, as a $0$--$1$ vector. By Farkas' lemma there is
$y\in\R^{\text{cells}}$ with $\langle y,\one_U\rangle\ge0$ for every such $U$
and $\langle y,\one\rangle<0$. Let $f$ be the $q$-step function with value
$y$ on each cell, and let $f'=f-c$ where $c=\int_Gf/|E|<0$. Then
$\int_G f'=0$ and $I_{f'}(U)=I_f(U)-cr\ge -cr>0$ for all $U\in\X_r^{(q)}(G)$,
hence for all $U\in\X_r(G)$ by Lemma~\ref{lem:grid}. Thus $r\notin\cH(G)$.

The last sentence follows from \eqref{eq:Prational}.
\end{proof}

\begin{corollary}
If $r=p/q\in\cP(G)$ in lowest terms, then $r$ is realized by an exact
uniform cover consisting of $q$-aligned sets.
\end{corollary}

\subsection{A duality characterization}
Let $K_r\subseteq L^1(G)$ be the convex cone generated by $\{\one_U:U\in\X_r(G)\}$.

\begin{theorem}\label{thm:dual}
For $0<r<|E|$, $r\in\cH(G)$ if and only if $\one\in\overline{K_r}$ (closure
in $L^1(G)$). Moreover, the class $\RI(G)$ in the definition of $\cH(G)$ may
be replaced by $L^\infty(G)$, or by the $q$-step functions for $q$ ranging over $\N$,
without changing $\cH(G)$.
\end{theorem}

\begin{proof}
Suppose $\one\in\overline{K_r}$ but $r\notin\cH(G)$. By
Lemma~\ref{lem:sign} there are $f\in\RI(G)$ with $\int_Gf=0$ and
$\varepsilon>0$ with $I_f\ge\varepsilon$ on $\X_r(G)$. For
$h=\sum\lambda_\alpha\one_{U_\alpha}\in K_r$ we get
$\int fh\ge\varepsilon\sum\lambda_\alpha=\varepsilon\int h/r$. Letting
$h\to\one$ in $L^1$ gives $0=\int f\ge\varepsilon|E|/r$, a contradiction.

Conversely, if $\one\notin\overline{K_r}$, the Hahn--Banach theorem
($L^1(G)^*=L^\infty(G)$) gives $\varphi\in L^\infty(G)$ with $\int\varphi h\ge0$
for $h\in K_r$ and $\int\varphi<0$. Put $c=-\int\varphi/|E|>0$ and
$\varphi'=\varphi+c$; then $\int\varphi'=0$ and $I_{\varphi'}(U)\ge cr$ on
$\X_r(G)$. Choose a step function $g$ with $\|g-\varphi'\|_{L^1}<cr/4$ and
let $g'=g-\int g/|E|$. Then $|\int g|<cr/4$, and for $U\in\X_r(G)$,
$I_{g'}(U)\ge cr-\tfrac{cr}{4}-\tfrac{cr}{4}>0$. So $r\notin\cH(G)$, witnessed by a step
function. The same argument shows that enlarging or shrinking the test class
between step functions and $L^\infty$ does not change $\cH(G)$.
\end{proof}

\begin{lemma}[Averaging]\label{lem:avg}
If $\nu$ is a Borel probability measure on $\X_r(G)$ such that
$\nu\{U:x\in U\}=r/|E|$ for almost every $x\in G$, then $r\in\cH(G)$.
\end{lemma}

\begin{proof}
For $\int_Gf=0$, Fubini gives $\int I_f\,d\nu=\int_Gf(x)\,\nu\{U\ni x\}\,dx=0$,
so $I_f$ cannot be positive on all of $\X_r(G)$; apply Lemma~\ref{lem:sign}.
\end{proof}

Combining \eqref{eq:Prational}, Theorem~\ref{thm:rational} and
Proposition~\ref{prop:closed}:

\begin{corollary}\label{cor:reduction}
$\overline{\cP(G)}=\overline{\cH(G)\cap\Q}$. Hence $\cH(G)=\overline{\cP(G)}$ if and only if
every irrational point of $\cH(G)$ is a limit of rational points of
$\cH(G)$; in particular it holds whenever $\cH(G)$ is a union of
nondegenerate intervals and rational points.
\end{corollary}

\begin{conjecture}\label{conj:PH}
For every finite connected graph, every irrational point of $\cH(G)$ is a
limit of rational points of $\cH(G)$; equivalently, $\cH(G)=\overline{\cP(G)}$.
\end{conjecture}

By Theorem~\ref{thm:dual} this asks whether an approximate uniform cover at
an irrational $r$ can always be converted into exact uniform covers at
rationals $r'\to r$.

\section{A covering lemma}\label{sec:covering}

The basic tool is the circle version of the sliding-window argument.

\begin{lemma}[Circle lemma]\label{lem:circle}
Let $\mathbb T_L=\R/L\Z$ and $h$ Riemann integrable on $\mathbb T_L$ with
$\int h=0$. For every $r\in[0,L]$ there is an arc $B$ of length $r$ with
$\int_B h=0$.
\end{lemma}

\begin{proof}
$F(x)=\int_x^{x+r}h$ is continuous on $\mathbb T_L$ and
$\int_0^LF(x)\,dx=r\int h=0$, so $F$ vanishes somewhere.
\end{proof}

A \emph{closed walk} $W$ of length $L$ in $G$ is a cyclic sequence of
$L$ edge traversals; it defines a map $\phi_W:\mathbb T_L\to G$ that is an
isometry on each unit interval $[j,j+1]$. The \emph{injectivity length}
$\iota(W)\in[0,L]$ is the largest $\rho$ such that $\phi_W$ is injective
outside a finite set on every arc of length $\rho$. A family
$W_1,\dots,W_m$ is a \emph{uniform cover} of multiplicity $k$ if every edge
is traversed exactly $k$ times in total.

\begin{lemma}[Covering lemma]\label{lem:covering}
If $G$ has a uniform cover by closed walks $W_1,\dots,W_m$ with
$\iota(W_i)\ge\rho$ for all $i$, then $[0,\rho]\subseteq\cH(G)$.
\end{lemma}

\begin{proof}
Let $0<r\le\rho$ and $\int_Gf=0$. Put $c_i=\int_{\mathbb T_{L_i}}f\circ\phi_{W_i}$.
Then $\sum_ic_i=k\int_Gf=0$, so $c_i\ge0\ge c_j$ for some $i,j$. Apply
Lemma~\ref{lem:circle} to $h=f\circ\phi_{W_i}-c_i/L_i$ on $\mathbb T_{L_i}$:
there is an arc $B$ of length $r$ with $\int_Bf\circ\phi_{W_i}=rc_i/L_i\ge0$.
Since $r\le\iota(W_i)$, $U=\phi_{W_i}(B)$ is a closed connected set of
measure $r$ and $I_f(U)=\int_Bf\circ\phi_{W_i}\ge0$. Similarly some
$U'\in\X_r(G)$ has $I_f(U')\le0$. Conclude with Lemma~\ref{lem:sign}.
\end{proof}

Note that the conclusion needs only one walk with $c_i\ge0$ and one with
$c_j\le0$; no relation between the lengths $L_i$ is required.

\begin{corollary}\label{cor:euler}
If $G$ is Eulerian, then $\cH(G)=[0,|E|]$.
\end{corollary}

\begin{proof}
An Euler circuit is a uniform cover of multiplicity $1$ by one walk with
$\iota=|E|$.
\end{proof}

A closed walk is \emph{non-backtracking} if no edge traversal is followed
(cyclically) by the reverse traversal of the same edge. A traversal of a
non-loop edge followed by any traversal of the same edge is necessarily a
reversal; for a loop, traversing it twice in the same direction is allowed.
An arc of length at most $1$ meets at most two consecutive traversals, so:

\begin{lemma}\label{lem:nb}
A non-backtracking closed walk $W$ has $\iota(W)\ge1$.
\end{lemma}

\begin{proof}
An arc of length $r\le1$ meets at most two consecutive traversals, in
pieces of lengths $a$ and $b$ with $a+b=r\le1$. If the traversals use
different edges, the images meet in at most a vertex. If they use the same
edge, it is a loop traversed twice in the same direction, and the images
are $[1-a,1]$ and $[0,b]$, which meet in at most one point.
\end{proof}

\begin{lemma}\label{lem:doublecover}
If $G$ is connected and $\deg v\ge2$ for all $v$, then $G$ has a uniform
cover of multiplicity $2$ by non-backtracking closed walks.
\end{lemma}

\begin{proof}
Call a subgraph $K$ \emph{good} if its edges admit such a cover. A cycle
(possibly a loop or a $2$-cycle) exists because $\deg\ge2$, and is good
(traverse it twice as two separate walks). Let $K$ be a maximal good
subgraph and $L=G-E(K)$. If $L$ contained a cycle, adding it twice would
contradict maximality, so $L$ is a forest. Suppose $L$ has an edge, and let
$\alpha$ be a longest path in $L$, from $u$ to $v$. Then $u\neq v$, and
$\deg_L u=\deg_L v=1$, so $\deg_G\ge2$ forces $u,v\in V(K)$. Let
$\gamma,\eta$ be walks of the cover through $u$ and $v$.
If $\gamma\neq\eta$, base $\gamma$ at $u$ and $\eta$ at $v$ and replace the
two walks by
\[
\gamma\cdot\alpha\cdot\eta\cdot\alpha^{-1}.
\]
If $\gamma=\eta$, write $\gamma=\gamma_1\gamma_2^{-1}$ with $\gamma_1,\gamma_2$
walks from $u$ to $v$, and replace $\gamma$ by
\[
\gamma_1\cdot\alpha^{-1}\cdot\gamma_2\cdot\alpha^{-1}
\]
At every junction a $K$-edge meets an $L$-edge, so no backtracking is created;
the edges of $\alpha$ are covered twice and all other counts are unchanged.
This contradicts maximality, so $L$ has no edges and $K=G$.
\end{proof}

\begin{theorem}\label{thm:mindeg}
If $G$ is connected and has no vertex of degree $1$, then $[0,1]\subseteq\cH(G)$.
\end{theorem}

\begin{proof}
Lemmas~\ref{lem:doublecover}, \ref{lem:nb} and \ref{lem:covering}.
\end{proof}

\begin{theorem}\label{thm:girth}
If $G$ is bridgeless and has girth $g$ (a loop has length $1$, a pair of
parallel edges forms a cycle of length $2$), then $[0,g]\subseteq\cH(G)$.
\end{theorem}

\begin{proof}
Every bridgeless graph has a uniform cover by cycles: by Seymour's
characterization of the cone generated by circuits \cite{Seymour}, the
all-ones vector is a nonnegative rational combination of circuits as soon
as no cut consists of a single edge; alternatively, Bermond, Jackson and
Jaeger \cite{BJJ} proved that every bridgeless graph has a cycle cover in
which each edge is covered exactly four times. A simple cycle of length
$n\ge g$ has $\iota=n\ge g$. Apply Lemma~\ref{lem:covering}.
\end{proof}

\begin{remark}\label{rem:unify}
Corollary~\ref{cor:euler}, Theorem~\ref{thm:mindeg} and Theorem~\ref{thm:girth}
are the cases $\rho=|E|$, $\rho=1$ and $\rho=g$ of Lemma~\ref{lem:covering}.
More generally, for an integer $\rho$, if every $\rho+1$ consecutive traversals of each $W_i$ use
distinct edges, then $\iota(W_i)\ge\rho$, and longer ``locally injective''
covers give longer initial intervals in $\cH(G)$.
\end{remark}

\section{Bridges}\label{sec:bridges}

\begin{proposition}\label{prop:bridge}
Let $e$ be a bridge of $G$ and let $m_1,m_2$ be the numbers of edges of
the two components of $G-e$, so $m_1+m_2=|E|-1$. Then
\[
\cH(G)\cap\bigl(\max\{|E|/2,\,m_1,\,m_2\},\ |E|\bigr)=\emptyset .
\]
\end{proposition}

\begin{proof}
Let $p\in e$ be at distance $t\in(0,1)$ from the $m_1$-side. Each component
of $G\setminus\{p\}$ has measure at most
$\mu_t=\max\{m_1+t,\,m_2+1-t\}$. Let $r\in(\mu_t,|E|)$ and choose
$\varepsilon>0$ with $\mu_t+\varepsilon<r$ and $J=[t-\varepsilon,t+\varepsilon]\subseteq e$.
If $U\in\X_r(G)$ missed a point $p'\in J$, it would lie in one component of
$G\setminus\{p'\}$, whose measure is at most $\mu_t+\varepsilon<r$. So every
$U\in\X_r(G)$ contains $J$. Let $f=\frac{c|E|}{2\varepsilon}\one_J-c$ with
$c>0$. Then $\int_Gf=0$ and $I_f(U)=c|E|-cr>0$ for every $U\in\X_r(G)$, so
$r\notin\cH(G)$. Finally,
$\inf_{0<t<1}\mu_t=\max\{|E|/2,m_1,m_2\}$.
\end{proof}

\begin{corollary}
If $\cH(G)=[0,|E|]$, then $G$ is bridgeless.
\end{corollary}

\begin{conjecture}\label{conj:bridgeless}
$\cH(G)=[0,|E|]$ if and only if $G$ is bridgeless.
\end{conjecture}

By Theorem~\ref{thm:rational} and Proposition~\ref{prop:closed}, the
``if'' direction is equivalent to the purely combinatorial statement that
for a bridgeless $G$ and every rational $r\in(0,|E|)$ there is an exact uniform
cover of $G$ by closed connected sets of measure $r$. Corollary~\ref{cor:euler}
and Theorem~\ref{thm:girth} prove it for Eulerian graphs and for
$r\le g(G)$. The next proposition shows how covers beyond the girth can be
built in a small case; the remaining evidence is computational.

\begin{proposition}[The theta graph]\label{prop:theta}
Let $\Theta_3$ have two vertices joined by three edges ($|E|=3$, girth
$2$, not Eulerian). Then $[0,2]\cup[\tfrac52,3]\subseteq\cH(\Theta_3)$, and every
rational $r=p/q$ with $q\le6$ lies in $\cH(\Theta_3)$.
\end{proposition}

\begin{proof}
$[0,2]$ comes from Theorem~\ref{thm:girth}. Let $r=3-s$ with
$0<s\le\frac12$. Closures of complements of the following open sets $W$ of
measure $s$ lie in $\X_r(\Theta_3)$: (i) an interval $[x,x+s]$ in one edge;
(ii) a star around a vertex with arms of lengths $t_1+t_2+t_3=s$. Choose a
random $W$ as follows. With probability $\pi_I$ take (i), with the edge
uniform and $x$ uniform on $[0,1-s]$. Otherwise choose one of the two
vertices uniformly and take (ii) with arm vector a uniformly random
permutation of $(0,t,s-t)$, $t$ uniform on $[0,s]$. Under (i), a point at
distance $y$ from the nearest vertex is covered with probability
$\frac13\min(y,s)/(1-s)$. Under (ii), each arm is $0$ with probability $\frac13$
and otherwise uniform on $[0,s]$, so the point is covered with probability
$\frac12\cdot\frac23(1-y/s)_+$. Taking $\pi_I/(1-\pi_I)=(1-s)/s$ makes the
total coverage probability independent of $y$, so the complements
$U=\overline{\Theta_3\setminus W}$ satisfy Lemma~\ref{lem:avg}. The last
assertion is the computation in Section~\ref{sec:comp}.
\end{proof}

\section{Computations}\label{sec:comp}

Theorem~\ref{thm:rational} turns $r\in\cH(G)$, for rational $r$, into a
finite LP. We enumerated $\X_r^{(q)}(G)$ and solved the LP for all
$r=p/q\in(0,|E|)$ with $q\le Q$. Each answer was certified in exact
rational arithmetic: a positive answer by a rational uniform cover, a
negative answer by a $q$-step function $f$ with $\int f=0$ and $I_f>0$ on
$\X_r^{(q)}(G)$ (hence on $\X_r(G)$ by Lemma~\ref{lem:grid}). In a few cases
the floating-point cover could not be rounded to an exact one; there the LP
was feasible to solver tolerance.

\begin{table}[h]
\centering\small
\begin{tabular}{@{}llcl@{}}
\toprule
graph & description & $Q$ & $r=p/q\in(0,|E|)$, $q\le Q$, \emph{not} in $\cH(G)$\\
\midrule
$I$ & one edge & 8 & all except $1/k$ \hfill (Lévy)\\
$P_2$ & path, 2 edges & 6 & all except $2/k$\\
$K_{1,3}$ & star & 6 & $\tfrac35,\tfrac45,\tfrac56$, and all in $(1,2)\cup(2,3)$\\
lollipop & loop $+$ pendant edge & 6 & all in $(1,2)$\\
dumbbell & two loops joined by a bridge & 6 & all in $(1,3)$ except $\tfrac32$\\
tadpole & triangle $+$ pendant edge & 6 &
  $\tfrac35,\tfrac76,\tfrac65,\tfrac75,\tfrac32,\tfrac74,\tfrac{11}6$, all in $(2,\tfrac52)$,\\
 & & & $\tfrac{11}4,\tfrac{14}5,\tfrac{17}6$, all in $(3,4)$\\
\midrule
$\Theta_3$ & 3 parallel edges & 6 & none\\
$\Theta_{1,1,2}$ & theta with one path of length 2 & 6 & none\\
$C_4+$chord & 4-cycle with a diagonal & 6 & none\\
$K_4^-$ & $K_4$ minus an edge & 6 & none\\
$K_4$ & complete graph on 4 vertices & 6 & none\\
two triangles & joined by two disjoint edges & 3 & none\\
\bottomrule
\end{tabular}
\caption{Exact chord set computations. The graphs below the rule are
bridgeless and non-Eulerian.}\label{tab:comp}
\end{table}

Observations.
\begin{enumerate}[label=(\arabic*)]
\item Every bridgeless example has every tested rational in $\cH(G)$,
  supporting Conjecture~\ref{conj:bridgeless}; the girth bound of
  Theorem~\ref{thm:girth} is far from sharp.
\item The gaps forced by Proposition~\ref{prop:bridge} are visible:
  $(\tfrac32,3)$ for the dumbbell, $(1,2)$ for the lollipop, $(2,3)$ for
  $K_{1,3}$, $(3,4)$ for the tadpole. For the lollipop the observed gap is
  exactly the predicted one; for the dumbbell $\frac32\in\cH$ shows that the
  bound is attained, but there are further gaps in $(1,\frac32)$.
\item The hypothesis of Theorem~\ref{thm:mindeg} is not necessary: the
  lollipop has a leaf, yet every tested $r\le1$ lies in its chord set.
  On the other hand $K_{1,3}$ and the tadpole miss $r=\frac35<1$.
\end{enumerate}

\section{Necklaces and a game}\label{sec:necklace}

Two thieves share an open necklace with $2N$ white and $2N$ black pearls;
two cuts always suffice to divide it fairly \cite{GW,Alon}. In the
continuous form (Totik \cite{Totik}): if $f,g\in\RI([0,1])$ with
$\int f=\int g=1$, then for every $r\in[0,1]$ there is an interval $J$ with
$\int_Jf=\int_Jg=r$ or $\int_Jf=\int_Jg=1-r$. Cutting the circle at a
point gives the circle version: for $f,g$ on $S^1$ with $\int f=\int g=1$
and every $r\in[0,1]$, some arc $U$ has $\int_Uf=\int_Ug=r$. With
$g\equiv1$ this is Lemma~\ref{lem:circle}, which is all that the results
of Sections~\ref{sec:covering}--\ref{sec:bridges} use.

\subsection*{The game}
On a circle with $2N$ dots, two players alternately cross out between $1$
and $m$ uncrossed dots, one in red and one in blue; once a player has
crossed out $N$ dots the other crosses out the rest. Whoever first creates
an arc of $2n$ consecutive dots with $n$ of each color loses. A loser
always exists, by the following discrete lemma: \emph{if $2N$ beads on a
circle are colored with $N$ red, then for every $1\le n\le N$ some $2n$
consecutive beads contain exactly $n$ red.} Indeed, the number $R_i$ of red
beads in the window starting at bead $i$ changes by at most $1$ when $i$
increases, and $\sum_iR_i=2nN$, so the average of $R_i$ is $n$.

\subsection*{The game on graphs}
The same argument works on a graph $G$ when the $2N$ dots are distinct points
in the interiors of edges and the losing configurations are closed connected
sets containing $2n$ dots, $n$ of each color: if $G$ is Eulerian, the Euler
circuit maps a circle onto $G$, bijectively away from the vertices, so arcs
become connected sets with the same dots.

If instead the dots are the \emph{vertices} of $G$, as the phrase ``connected
subset containing $2n$ dots'' suggests, the Euler-circuit argument fails
because an Euler circuit revisits vertices. Indeed, let $G$ be the path
$v_1\cdots v_8$ with every edge doubled (an Eulerian multigraph) and color
the vertices $BBRRRRBB$. The connected $6$-sets are the three windows of
length $6$, each with four red vertices, so with $N=4$, $n=3$ the game ends
without a loser. An exhaustive search shows that every connected graph with
at most $6$ vertices always has a loser in the vertex version, so $8$
vertices is the minimum. A Hamiltonian cycle suffices for the vertex version.

\section{Further questions}

\begin{enumerate}[label=(\arabic*)]
\item Prove Conjecture~\ref{conj:bridgeless}. A natural route is to build,
  for a bridgeless $G$ and each $r$, a random closed connected set of
  measure $r$ whose coverage probability is the same at every point
  (Lemma~\ref{lem:avg}); Proposition~\ref{prop:theta} does this for $\Theta_3$
  and $r\ge\frac52$.
\item Prove or disprove Conjecture~\ref{conj:PH}. Is the complement of
  $\cH(G)$ in $(0,|E|)$ always a finite union of open intervals together
  with countably many points accumulating only at rational numbers, as for
  the interval?
\item Find a Hopf-type structure theorem: for the interval, the complement
  of the chord set in $(0,1)$ is closed under addition. Is there an
  analogue for graphs, perhaps relative to the edge lengths of cycles?
\item Describe $\cH(G)$ for trees. For a path of length $n$ it is
  $\{n/k\}\cup\{0\}$. The star $K_{1,3}$ is already irregular
  (Table~\ref{tab:comp}).
\item Metric graphs with arbitrary edge lengths. All definitions and the
  results of Sections~\ref{sec:prelim}--\ref{sec:bridges} extend
  verbatim (the grid lemma for commensurable lengths). Does $\cH(G)$ vary
  upper semicontinuously with the edge lengths?
\item Which graphs satisfy $[0,1]\subseteq\cH(G)$? The lollipop shows that
  minimum degree $2$ is not necessary.
\end{enumerate}

\section*{Acknowledgement}
This revised version is prepared by an AI assistant (Claude).


\begin{thebibliography}{99}
\bibitem{Alon} N.~Alon, Splitting necklaces, \emph{Adv. Math.} 63 (1987), 247--253.
\bibitem{BJJ} J.-C.~Bermond, B.~Jackson, F.~Jaeger, Shortest coverings of graphs with cycles, \emph{J. Combin. Theory Ser. B} 35 (1983), 297--308.
\bibitem{GW} C.~H.~Goldberg, D.~B.~West, Bisection of circle colorings, \emph{SIAM J. Algebraic Discrete Methods} 6 (1985), 93--106.
\bibitem{Hopf} H.~Hopf, Über die Sehnen ebener Kontinuen und die Schleifen geschlossener Wege, \emph{Comment. Math. Helv.} 9 (1937), 303--319.
\bibitem{Levy} P.~Lévy, Sur une généralisation du théorème de Rolle, \emph{C. R. Acad. Sci. Paris} 198 (1934), 424--425.
\bibitem{Oxtoby} J.~C.~Oxtoby, Horizontal chord theorems, \emph{Amer. Math. Monthly} 79 (1972), 468--475.
\bibitem{Schrijver} A.~Schrijver, \emph{Theory of Linear and Integer Programming}, Wiley, 1986.
\bibitem{Seymour} P.~D.~Seymour, Sums of circuits, in \emph{Graph Theory and Related Topics} (J.~A.~Bondy, U.~S.~R.~Murty, eds.), Academic Press, 1979, 341--355.
\bibitem{Totik} V.~Totik, A tale of two integrals, \emph{Amer. Math. Monthly} 106 (1999), 227--240.
\end{thebibliography}
\end{document}